\documentclass[12pt]{amsart}
\usepackage{amsmath, amssymb, latexsym, amsthm, mathrsfs, color, mathtools, tikz,pgfplots,tikz-cd, graphicx, stackengine, url}
\usepackage[top=2.5cm, bottom=2.5cm, left=2.5cm, right=2.5cm]{geometry}

\stackMath
\newcommand\tsup[2][2]{%
 \def\useanchorwidth{T}%
  \ifnum#1>1%
    \stackon[-1pt]{\tsup[\numexpr#1-1\relax]{#2}}{\hspace{1pt}\scriptstyle\sim}%
  \else%
    \stackon[.5pt]{#2}{\hspace{1pt}\scriptstyle\sim}%
  \fi%
}

\newcommand{\nc}{\newcommand}

\newcommand{\swg}{\mathsf{S}_1(\Omega,\Gamma)}
\newcommand{\sgg}{{\mathsf{S}_1(\Ga,\Ga)}}
\newcommand{\sbgg}{{\mathsf{S}_1(\Ga_\mathrm{Bor},\Ga_\mathrm{Bor})}}

\newcommand{\sww}{\mathsf{S}_1(\Omega,\Omega)}
\newcommand{\swo}{\mathsf{S}_1(\Omega,\Op)}
\newcommand{\soo}{\mathsf{S}_1(\Op,\Op)}
\DeclareMathOperator{\C}{C_p}

\nc{\mc}{\mathcal}
\nc{\thusfar}{\my{--- Edited thus far ---}}
\nc{\lei}{\le^\oo}
\nc{\sqsubs}{\sqsubseteq^*}
\nc{\card}[1]{\left|#1\right|}
\nc{\medcard}[1]{\biggl|\,#1\,\biggr|}
\nc{\smallcard}[1]{|\,#1\,|}
\nc{\bds}{bidirectional $\roth$-scale}
\nc{\bfP}{\mathbf{P}}
\nc{\bfQ}{\mathbf{Q}}
\nc{\bbT}{\mathbb{T}}
\nc{\bbZ}{\mathbb{Z}}
\nc{\bbN}{\mathbb{N}}
\nc{\bbC}{\mathbb{C}}
\nc{\beq}{\begin{equation}}\nc{\eeq}{\end{equation}}
\nc{\mbq}{\mb{?}}
\nc{\mb}[1]{{\mbox{\textbf{#1}}}}
\nc{\nop}{$\times$}
\nc{\fbn}{\!\!\fbox{\!\nop\!}\!\!}
\nc{\yup}{\checkmark}
\nc{\forces}{\Vdash}
\nc{\name}[1]{\dot{#1}}
\nc{\tf}{\my{FINISHED THUS FAR}}
\nc{\FU}{Fr\'echet--Urysohn}
\nc{\gs}{$\gamma$~space}
\nc{\Ga}{\Gamma}
\nc{\Gab}{\Ga_\mathrm{Bor}}
\nc{\Om}{\Omega}

\nc{\smallbinom}[2]{\begin{psmallmatrix} #1\\ #2 \end{psmallmatrix}}
\nc{\bgamma}{\smallbinom{\Om}{\Ga}}
\newcommand{\two}{\{0,1\}}
\nc{\productive}[2]{(#1,\allowbreak #2)^\x}
\nc{\prdct}[1]{#1^\x}
\nc{\Sel}{\mathsf{S}}
\nc{\sset}[2]{\{\,#1 : #2\,\}}
\nc{\smb}[1]{{\!\!\mb{#1}\!\!}}
\nc{\medset}[2]{{\biggl\{\,#1 : #2\,\biggr\}}}
\nc{\smallmedset}[2]{{\bigl\{\,#1 : #2\,\bigr\}}}
\nc{\set}[2]{{\left\{\,#1 : #2\,\right\}}}
\nc{\eseq}[1]{#1_1, \allowbreak #1_2, \allowbreak\dotsc} 
\nc{\seleseq}[1]{#1_1\in \mathcal{#1}_1, \allowbreak #1_2\in \mathcal{#1}_2, \allowbreak\dotsc}
\nc{\cube}{(\Cantor)^\bbN}
\nc{\Match}{\op{Match}}
\nc{\concat}[1]{\hat{\phantom{a}}\langle #1\rangle}
\nc{\poset}{\mathbb{P}}
\nc{\fn}[1]{{\op{Fn}(#1\times\w,2)}}
\nc{\linadd}{\op{linadd}}
\nc{\nonprod}{\non^\x}
\nc{\alephes}{{\aleph_0}}
\nc{\my}[1]{\marginpar{\textcolor{red}{***}}\textcolor{red}{#1}}
\nc{\myb}[1]{\marginpar{\textcolor{blue}{***}}\textcolor{blue}{#1}}
\nc{\later}[1]{{\color{green} #1}}
\nc{\BTs}[1]{{\color{green} #1 (BT)}}
\nc{\Cp}{\op{C}_\mathrm{p}}
\nc{\Bp}{\op{B}_p}
\nc{\Pa}[8]{\bibitem{#1} {#2}, \emph{#3}, {#4} \textbf{#5} ({#6}), {#7}--{#8}.}
\nc{\tPa}[5]{\bibitem{#1} {#2}, \emph{#3}, {#4}, to appear.}
\nc{\sPa}[4]{\bibitem{#1} {#2}, \emph{#3}, {#4}, submitted.}
\nc{\Bc}[9]{\bibitem{#1} {#2}, \emph{#3}, in: \textbf{#4} (#5), #6 #7, #8--#9.}
\nc{\fD}{\mathfrak{D}}
\nc{\fX}{\mathfrak{X}}
\nc{\Onbd}{\Op_{\mathrm{nbd}}} 
\nc{\Omnb}{\Om_{\mathrm{nbd}}} 
\nc{\od}{\mathfrak{od}}
\nc{\Setting}[7]{\xymatrix@R=4pt@C=7pt{#1\ar@{-}[r]&#2\ar@{-}[r]&#3\\&#4\ar@{-}[u]\\
#5\ar@{-}[uu]\ar@{-}[r] & #6\ar@{-}[u]\ar@{-}[r] & #7\ar@{-}[uu]}}
\nc{\mx}[1]{\begin{matrix}#1\end{matrix}}
\nc{\plim}{p\txt{-}\lim}
\nc{\Bgp}{{\Z^\bbN}}
\nc{\Cgp}{{{\Z_2}^\bbN}}
\nc{\Cite}[1]{\textbf{[#1]}}
\nc{\Next}[1]{{#1^+}}
\nc{\cFin}{\mathrm{cF}}
\nc{\scsp}{\text{-scale space}}
\nc{\cfn}{\text{cofinal}\ }
\nc{\Con}{\text{Concentrated}}
\nc{\Lind}{\text{Lindel\"of}\,}
\nc{\con}{\text{-Concentrated}}
\nc{\lind}{\text{-Lindel\"of}\,}
\nc{\ctbl}{\text{countably }\allowbreak}
\nc{\Hur}{\text{Hurewicz}}
\nc{\intvl}[2]{{[#1(#2),\allowbreak #1(#2\!+\!1))}}
\nc{\Bdd}{\mathbf{B}}
\nc{\Dfin}{\mathfrak{D}_\mathrm{fin}}
\nc{\grbl}{{\mbox{\textit{\tiny gp}}}}
\nc{\bbP}{\mathbb{P}}
\nc{\BOfat}{\B_{\Om_{\mathrm{fat}}}}
\nc{\Bgood}{\B_{\mathrm{good}}}
\nc{\compactN}{\cl{\mathbb{N}}}
\nc{\blocks}[2]{\op{cl}_{#2}(#1)}
\nc{\blocksplus}[2]{\op{cl}^+_{#2}(#1)}
\nc{\arx}[1]{\texttt{http://arxiv.org/math/#1}}
\nc{\bq}{\begin{quote}}
\nc{\eq}{\end{quote}}
\nc{\cl}[1]{\overline{#1}}
\nc{\Cl}[2]{\mathrm{cl}_{#1}(#2)}
\nc{\CH}{the Continuum Hypothesis}
\nc{\MA}{Martin's Axiom}
\nc{\Bfat}{\B_\mathrm{fat}}
\nc{\inv}{^{-1}}
\nc{\Cantor}{{\two^\bbN}}
\nc{\bP}{\mathbf{P}}
\nc{\bof}{\op{\fb}}
\nc{\dof}{\op{\fd}}
\nc{\bofF}{\bof(\cF)}
\nc{\sr}[3]{\underset{\mbox{#3}}{\mbox{#1}}}
\nc{\gp}{\binom{\Om}{\Ga}}
\nc{\gpsmall}{\mbox{$\gp$}}
\nc{\gig}{\gimel}
\nc{\gns}{\sone(\Om,\gig)}
\nc{\nsr}[2]{#1}
\nc{\Srg}{{\mathbb{S}}}
\nc{\Srgs}{{\mathbb{S}^*}}
\nc{\NN}{{\bbN^{\bbN}}}
\nc{\ZN}{{\Z^{\bbN}}}
\nc{\NNup}{{\bbN^{\uparrow\bbN}}}
\nc{\NNupb}{{b^{\uparrow\bbN}}}
\nc{\Pof}{\op{P}}
\nc{\PN}{{\Pof(\bbN)}}
\nc{\rothel}[1]{{[#1]^{\mbox{\tiny $\infty$}}}}
\nc{\roth}{{[\bbN]^{\mbox{\tiny $\infty$}}}} 
\nc{\roths}{{[b]^{\mbox{\tiny $\infty$}}}} 
\nc{\Fin}{\mathrm{Fin}}
\nc{\ici}{[\bbN]^{ \infty, \infty}}
\nc{\Inc}{{\compactN^{\uparrow\bbN}}}
\nc{\powInc}[1]{{\big(\Inc\big)^{#1}}}
\nc{\powFin}[1]{{\big(\Fin\big)^{#1}}}
\nc{\powPN}[1]{{\big(\PN\big)^{#1}}}
\nc{\NcompactN}{{\compactN^\bbN}}
\nc{\seq}[2]{\la #1\ra_{#2\in\bbN}}
\nc{\Uarrow}{\smash{\big\uparrow}}
\nc{\LE}{\preccurlyeq}
\nc{\GE}{\succcurlyeq}
\nc{\op}{\operatorname}
\nc{\im}{\op{im}}
\nc{\Span}{\op{span}}
\nc{\maxfin}{\op{maxfin}}
\nc{\ran}{\op{range}}
\nc{\iso}{\cong}
\nc{\Madd}{{\M}^\star}
\nc{\cI}{\mathcal{I}}
\nc{\cJ}{\mathcal{J}}
\nc{\scrA}{\mathscr{A}}
\nc{\scrB}{\mathscr{B}}
\nc{\scrC}{\mathscr{C}}
\nc{\scrD}{\mathscr{D}}
\nc{\scrF}{\mathscr{F}}
\nc{\scrK}{\mathscr{K}}
\nc{\A}{\D\forall}
\nc{\B}{\mathrm{B}}
\nc{\cB}{\mathcal{B}}
\nc{\cZ}{\mathcal{Z}}
\nc{\bB}{\mathbf{B}}
\nc{\BS}{\mathbf{B}(\mathcal{S})}
\nc{\BF}{\mathbf{B}(\mathcal{F})}
\nc{\BU}{\mathbf{B}(\mathcal{U})}
\nc{\cSp}{\mathcal{S}^+}
\nc{\cFp}{\mathcal{F}^+}
\nc{\cUp}{\mathcal{U}^+}

\nc{\BG}{\B_\Ga}
\nc{\BL}{\B_\Lambda}
\nc{\BT}{\B_\Tau}
\nc{\BTstar}{\B_{\Tau^*}}
\nc{\BO}{\B_\Om}
\nc{\DO}{\cD_\Om}
\nc{\KO}{\cK_\Om}
\nc{\CG}{C_\Ga}
\nc{\CL}{C_\Lambda}
\nc{\CT}{C_\Tau}
\nc{\CTstar}{C_{\Tau^*}}
\nc{\CO}{C_\Om}
\nc{\COgp}{C_{\Om^{\grbl}}}
\nc{\CLgp}{C_{\Lambda^{\grbl}}}
\nc{\BOgp}{\B_{\Om}^{\grbl}}
\nc{\BLgp}{\B_{\Lambda^{\grbl}}}
\nc{\sfC}{\mathsf{C}}
\nc{\sfD}{\mathsf{D}}
\nc{\bD}{\mathbf{D}}
\nc{\Tau}{\mathrm{T}}
\nc{\cA}{\mathcal{A}}
\nc{\cK}{\mathcal{K}}
\nc{\cD}{\mathcal{D}}
\nc{\cF}{\mathcal{F}}
\nc{\cS}{\mathcal{S}}
\nc{\cT}{\mathcal{T}}
\nc{\cG}{\mathcal{G}}
\nc{\cY}{\mathcal{Y}}
\nc{\J}{\mathcal{J}}
\nc{\cL}{\mathcal{L}}
\nc{\cM}{\mathcal{M}}
\nc{\cN}{\mathcal{N}}
\nc{\cH}{\mathcal{H}}
\nc{\cO}{\mathcal{O}}
\nc{\Op}{\mathrm{O}}
\nc{\rmA}{\mathrm{A}}
\nc{\rmF}{\mathrm{F}}
\nc{\rmB}{\mathrm{B}}
\nc{\rmD}{\mathrm{D}}
\nc{\rmP}{\mathrm{P}}
\nc{\cC}{\mathcal{C}}
\nc{\cP}{\mathcal{P}}
\nc{\bbQ}{\mathbb{Q}}
\nc{\bbR}{\mathbb{R}}
\nc{\cU}{\mathcal{U}}
\nc{\cQ}{\mathcal{Q}}
\nc{\Un}{\bigcup}
\nc{\cV}{\mathcal{V}}
\nc{\cW}{\mathcal{W}}
\nc{\Z}{{\mathbb Z}}
\nc{\Impl}{\Rightarrow}
\long\def\forget#1\forgotten{\marginpar{\textcolor{green}{Forgetting...}}}
\nc{\ft}{\mathfrak{t}}
\nc{\fb}{\mathfrak{b}}
\nc{\fc}{\mathfrak{c}}
\nc{\fd}{\mathfrak{d}}
\nc{\fg}{\mathfrak{g}}
\nc{\oo}{\infty}
\nc{\fr}{\mathfrak{r}}
\nc{\fk}{\mathfrak{k}}
\nc{\bidi}{\mathfrak{bidi}}
\nc{\fu}{\mathfrak{u}}
\nc{\fh}{\mathfrak{h}}
\nc{\fp}{\mathfrak{p}}
\nc{\fj}{\mathfrak{j}}
\nc{\fs}{\mathfrak{s}}
\nc{\w}{\omega}
\nc{\x}{\times}
\nc{\Iff}{\Leftrightarrow}
\newcommand\comp{^{\text{\tt c}}}
\nc{\nin}{\notin}
\nc{\cat}{\hat{\ }}
\nc{\sub}{\subseteq}
\nc{\spst}{\supseteq}
\nc{\sm}{\setminus}
\nc{\as}{\subseteq^*}
\nc{\les}{\le^*}
\nc{\leinf}{\le^{\infty}}
\nc{\leS}{\le_S}
\nc{\leF}{\le_F}
\nc{\leU}{\le_U}
\nc{\rest}{\restriction}

\nc{\la}{\langle}
\nc{\ra}{\rangle}
\nc{\E}{\exists}
\nc{\dom}{\op{dom}}
\nc{\cov}{\op{cov}}
\nc{\add}{\op{add}}
\nc{\addmen}{\add(\Men{})}
\nc{\cof}{\op{cof}}
\nc{\cf}{\op{cf}}
\nc{\non}{\op{non}}
\nc{\unif}{\op{non}}
\nc{\COV}{\op{COV}}
\nc{\ADD}{\op{ADD}}
\nc{\COF}{\op{COF}}
\nc{\NON}{\op{NON}}
\nc{\impl}{\to}
\nc{\Lp}{\mathcal{L_\p}}
\nc{\Wlog}{without loss of generality}
\newtheorem{thm}{Theorem}[section]
\nc{\bthm}{\begin{thm}} \nc{\ethm}{\end{thm}}
\newtheorem{prop}[thm]{Proposition}
\nc{\bprp}{\begin{prop}} \nc{\eprp}{\end{prop}}
\newtheorem{fact}[thm]{Fact}
\nc{\bfct}{\begin{fact}} \nc{\efct}{\end{fact}}
\newtheorem{prob}[thm]{Problem}
\nc{\bprb}{\begin{prob}} \nc{\eprb}{\end{prob}}
\newtheorem{lem}[thm]{Lemma}
\nc{\blem}{\begin{lem}} \nc{\elem}{\end{lem}}
\newtheorem{app}[thm]{Application}
\nc{\bapp}{\begin{app}} \nc{\eapp}{\end{app}}
\newtheorem{claim}[thm]{Claim}
\nc{\bclm}{\begin{claim}} \nc{\eclm}{\end{claim}}
\newtheorem{cor}[thm]{Corollary}
\nc{\bcor}{\begin{cor}} \nc{\ecor}{\end{cor}}
\newtheorem{conj}[thm]{Conjecture}
\nc{\bcnj}{\begin{conj}} \nc{\ecnj}{\end{conj}}
\theoremstyle{definition}
\newtheorem{defn}[thm]{Definition}
\nc{\bdfn}{\begin{defn}} \nc{\edfn}{\end{defn}}
\newtheorem{obs}[thm]{Observation}
\nc{\bobs}{\begin{obs}} \nc{\eobs}{\end{obs}}
\theoremstyle{remark}
\newtheorem{rem}[thm]{Remark}
\nc{\brem}{\begin{rem}} \nc{\erem}{\end{rem}}
\newtheorem{cnv}[thm]{Convention}
\nc{\bcnv}{\begin{cnv}} \nc{\ecnv}{\end{cnv}}
\newtheorem{exam}[thm]{Example}
\nc{\bexm}{\begin{exam}} \nc{\eexm}{\end{exam}}
\nc{\bpf}{\begin{proof}} \nc{\epf}{\end{proof}}
\nc{\be}{\begin{enumerate}}
\nc{\ee}{\end{enumerate}}
\nc{\bi}{\begin{itemize}}
\nc{\bimy}{\my{\begin{itemize}}
\nc{\eimy}{\end{itemize}}}
\nc{\itm}{\item}
\nc{\ei}{\end{itemize}}
\nc{\Subsection}[1]{\goodbreak\subsection*{#1}}

\nc{\sone}{\mathsf{S}_1}
\nc{\sfin}{\mathsf{S}_\mathrm{fin}}
\nc{\ufin}{\mathsf{U}_\mathrm{fin}}
\nc{\Split}{\mathsf{Split}}

\nc{\gone}{\mathsf{G}_1}    \nc{\gfin}{\mathsf{G}_\mathrm{fin}}

\nc{\men}{\sfin(\Op,\Op)}
\nc{\sch}{\ufin(\Op,\Omega)}
\nc{\rothb}{\sone(\Op,\Op)}
\nc{\pmen}{\sfin(\Omega,\Omega)}
\nc{\Rothb}{\sone(\Op,\Op)}

\nc{\prothb}{\sone(\Omega,\Omega)}
\nc{\tU}{{\tilde{U}}}
\nc{\tF}{{\tilde{F}}}
\nc{\tY}{{\tilde{Y}}}
\nc{\tX}{S}
\nc{\dtX}{X\sm S}
\nc{\dt}[1]{{\tsup[2]{#1}}}
\nc{\td}{{\tilde{d}}}
\nc{\tb}{{\tilde{b}}}
\nc{\tz}{{\tilde{z}}}
\nc{\cfd}{\cf(\fd)}
\nc{\msep}{\sfin(\cD,\cD)}
\nc{\rsep}{\sone(\cD,\cD)}
\nc{\cft}{\sfin(\Omega_{\mathbf{0}},\Omega_{\mathbf{0}})}
\nc{\scft}{\sone(\Omega_{\mathbf{0}},\Omega_{\mathbf{0}})}

\nc{\Umen}{U\text{-Menger}}
\nc{\hur}{\ufin(\cO,\Gamma)}
\nc{\tUmen}{\tU\text{-Menger}}

\nc{\Men}{\text{Menger}}
\nc{\Sch}{\text{Scheepers}}

\nc{\aspst}{\prescript{*}{}{\spst}\ }
\nc{\eqs}{=^*}
\nc{\ctblOm}{\Omega_{\mathrm{ctbl}}}
\nc{\GNga}{{\smallbinom{\Om}{\Ga}}}
\nc{\ctblga}{\smallbinom{\ctblOm}{\Ga}}
\nc{\Gaclp}{\Ga_{\mathrm{clp}}}
\nc{\sclpgg}{\sone(\Gaclp,\Gaclp)}

\nc{\sep}{
\vspace{2cm}
\noindent
\begin{minipage}{\textwidth}
	\textcolor{red}{\rule{\textwidth}{1pt}}
\end{minipage}
}

\nc{\bfzero}{\mathbf{0}}
\usepackage[all]{xy}

\title{Unbounded towers and products}

\author[P.~Szewczak]{Piotr Szewczak}
\address{Piotr Szewczak, Institute of Mathematics, Faculty of Mathematics and Natural Science College of Sciences, Cardinal Stefan Wyszy\'nski University in Warsaw, W\'oycickiego 1$\slash$3, 01--938 Warsaw, Poland
}
\email{p.szewczak@wp.pl}
\urladdr{http://piotrszewczak.pl}

\author[M.~W\l udecka]{Magdalena W\l udecka}
\address{Magdalena W\l udecka, Institute of Mathematics, Faculty of Mathematics and Natural Science College of Sciences, Cardinal Stefan Wyszy\'nski University in Warsaw, W\'oycickiego 1$\slash$3, 01--938 Warsaw, Poland}
\email{m.wludecka@gmail.com}

\keywords{unbounded tower, $\sgg$, Gerlits--Nagy, $\gamma$-property, $\gamma$-set, selection principles, products.}

\subjclass[2010]{26A03, 
54D20, 03E75, 03E17.
}

\begin{document}

\maketitle

\begin{abstract}
We investigate products of sets of reals with combinatorial covering properties.
A topological space satisfies $\sgg$ if for each sequence of point-cofinite open covers of the space, one can pick one element from each cover and obtain a point-cofinite cover of the space.
We prove that, if there is an unbounded tower, then there is a nontrivial set of reals satisfying $\sgg$ in all finite powers.
In contrast to earlier results, our proof does not require any additional set-theoretic assumptions.

A topological space satisfies $\GNga$ (also known as Gerlits--Nagy's property $\gamma$) if every open cover of the space such that each finite subset of the space is contained in a member of the cover, contains a point-cofinite cover of the space.
We investigate products of sets satisfying $\GNga$ and their relations to other classic combinatorial covering properties.
We show that finite products of sets with a certain combinatorial structure satisfy $\GNga$ and give necessary and sufficient conditions when these sets are productively $\GNga$.
\end{abstract}

\section{Introduction}

Let $\bbN$ be the set of natural numbers and $\roth$ be the set of infinite subsets of $\bbN$.
In a natural way an element of $\roth$ can be viewed as an increasing function in $\NN$.
A subset of $\roth$ is \emph{unbounded} if for each element $a\in\roth$, there is an element $b$ in the set such that the set $\sset{n}{a(n)>b(n)}$ is infinite.
A set $a$ is an \emph{almost subset} of a set $b$, denoted $a\as b$, if the set $a\sm b$ is finite.

\bdfn
Let $\kappa$ be an uncountable ordinal number.
A set $\sset{x_\alpha}{\alpha<\kappa}\sub\roth$ is a \emph{$\kappa$-unbounded tower} if it is unbounded and, for all ordinal numbers $\alpha,\beta<\kappa$ with $\alpha<\beta$, we have $x_\alpha\aspst x_\beta$.
\edfn

For an uncountable ordinal number $\kappa$, the existence of a $\kappa$-unbounded tower is independent of ZFC~\cite{BBC}.
It turns out that such a set is a significant object used in constructions nontrivial sets with combinatorial covering properties.

By \emph{space} we mean a topological space.
A \emph{cover} of a space is a family of proper subsets of the space whose union is the entire space, and an \emph{open} cover is a cover whose members are open subsets of the space.
A \emph{$\gamma$-cover} of a space is an infinite cover such that each point of the space belongs to all but finitely many sets of the cover.
For a space $X$, let $\Ga(X)$ be the family of all open $\gamma$-covers of the space $X$.
A space $X$ satisfies $\sgg$ if for each sequence $\eseq{\cU}\in\Ga(X)$, there are sets $\seleseq{U}$ such that $\sset{U_n}{n\in\bbN}\in\Ga(X)$.
This property was introduced by Scheepers~\cite{coc1} and it was studied in the context of local properties of functions spaces~\cite{SchCp, alpha_i}.

A \emph{set of reals} is a space homeomorphic to a subset of the real line.
We restrict our consideration to the realm of sets of reals.
We identify the family $\PN$ of all subsets of the set $\bbN$ with the Cantor space $\Cantor$, via characteristic functions.
Since the Cantor space is homeomorphic to Cantor's set, every subspace of the space $\PN$ is considered as a set of reals.
Let $\Fin$ be the family of all finite subsets of $\bbN$.
Let $\kappa$ be an uncountable ordinal number.
A set $X\cup\Fin$ is a \emph{$\kappa$-unbounded tower set}, if the set $X$ is a $\kappa$-unbounded tower in $\roth$.
Let $\fb$ be the minimal cardinality of an unbounded subset of $\roth$
Each $\fb$-unbounded tower set satisfies $\sone(\Ga,\Ga)$
(\cite[Theorem 6]{alpha_i}, \cite[Proposition 2.5]{BBC}).
Miller and Tsaban proved that, assuming some additional set-theoretic assumption, each $\fb$-unbounded tower set satisfies $\sgg$ in all finite powers~\cite[Theorem~2.8]{BBC}.
We prove that the same statement remains true with no extra assumption used by Miller and Tsaban. 
The proof method is new.

An \emph{$\w$-cover} of a space is an open cover such that, each finite subset of the space is contained in a set from the cover. 
A space satisfies $\GNga$ if every $\w$-cover of the space contains a $\gamma$-cover of the space.
This property was introduced by Gerlits and Nagy~\cite{gn}.
A \emph{pseudointersection} of a family of infinite sets is an infinite set $a$ with $a\as b$ for all sets $b$ in the family. A family of infinite sets is \emph{centered} if the finite intersections of its elements, are infinite. Let $\fp$ be the minimal cardinality of a subfamily of $\roth$ that is centered and has no pseudointersection.
By the result of Orenshtein and Tsaban (\cite[Theorem~3.6]{ot},~\cite[Theorem~6]{sss}), each $\fp$-unbounded tower set satisfies $\GNga$.
A set is \emph{productively $\GNga$} if its product space with any set satisfying $\GNga$, satisfies $\GNga$.
Miller, Tsaban and Zdomskyy proved that, each $\w_1$-unbounded tower set is productively $\GNga$~\cite[Theorem~2.8]{gamma}.
We show that, if each set of cardinality less than $\fp$ is productively $\GNga$, then each $\fp$-unbounded tower set is productively $\GNga$.
Moreover, the product space of finitely many $\fp$-unbounded tower sets, satisfies $\GNga$.
We also consider products of sets satisfying $\GNga$ and their relations to classic combinatorial covering properties.

\section{Generalized towers}

We generalize the notion of $\kappa$-unbounded tower for an uncountable ordinal number $\kappa$.
Let $n,m$ be natural numbers with $n<m$.
Define $[n,m):=\sset{i\in\bbN}{n\leq i<m}$.
A set $x\in \roth$ \emph{omits} the interval $[n,m)$, if $x\cap [n,m)=\emptyset$.

\blem[{Folklore~\cite[Lemma~2.13]{MHP}}]\label{lem:unbdd}
A set $X\sub\roth$ is unbounded if and only if for each function $a\in\roth$, there is a set $x\in X$ that omits infinitely many intervals $[a(n),a(n+1))$.\qed
\elem

Let $\kappa$ be an uncountable ordinal number, and $\set{x_\alpha}{\alpha<\kappa}$ be a $\kappa$-unbounded tower in $\roth$.
Fix a function $a\in\roth$.
By Lemma~\ref{lem:unbdd}, there is an ordinal number $\alpha<\kappa$ such that the set $x_\alpha$ omits infinitely many intervals $[a(n),a(n+1))$.
Thus, the set 
\[
b:=\sset{n\in\bbN}{x_\alpha\cap [a(n),a(n+1))=\emptyset}
\]
is an element of $\roth$ and
\[
x_\alpha\cap\Un_{n\in b}[a(n),a(n+1))=\emptyset.
\]
For each ordinal number $\beta<\kappa$ with $\alpha<\beta$, we have $x_\alpha\aspst x_\beta$, and thus 
\[
x_\beta\cap\Un_{n\in b}[a(n),a(n+1))\in\Fin.
\]
This observation motivates the following definition.

\bdfn
Let $\kappa$ be an uncountable ordinal number.
A set $X\sub \roth$ with $\card{X}\geq \kappa$ is a \emph{$\kappa$-generalized tower} if for each function $a\in\roth$, there are sets $b\in\roth$ and $S\sub X$ with $\card{S}<\kappa$ such that 
\[
x\cap \Un_{n\in b}[a(n),a(n+1))\in\Fin
\]
for all sets $x\in X\sm S$.
\edfn

Let $\kappa$ be an uncountable ordinal number.
Every $\kappa$-unbounded tower in $\roth$ is a $\kappa$-generalized tower.
The forthcoming Lemma~\ref{lem:uscale} shows that the notion of $\kappa$-unbounded tower may capture wider class of sets than $\kappa$-unbounded towers.

A set $B\sub\roth$ is \emph{groupwise dense} if:
\be
\item for all sets $y\in\roth$ and $b\in B$, if $y\as b$, then $y\in B$,
\item for each function $b\in\roth$, there is a set $c\in\roth$ such that $\Un_{n\in c}[b(n),b(n+1))\in B$.
\ee
The \emph{groupwise density number} $\fg$ is the minimal cardinality of a family of groupwise dense sets in $\roth$ with empty intersection. For function $f,g\in\NN$, let $f\circ g\in\NN$ be a function such that $(f\circ g)(n):=f(g(n))$ for all natural numbers $n$.

\blem\label{lem:uscale}
Let $\kappa$ be an uncountable ordinal number and $\lambda$ be an ordinal number with $\lambda<\fg$.
Let $X:=\Un_{\alpha<\lambda} X_\alpha$ be a union of $\kappa$-generalized towers.
\be
\item
If $\lambda<\cf(\kappa)$, then $X$ is a $\kappa$-generalized tower.
\item
If $\lambda\geq\cf(\kappa)$, then $X$ is a $(\kappa\cdot\lambda)^+$-generalized tower.
\ee
\elem

\bpf
(1) 
Fix a function $a\in\roth$ and an ordinal number $\alpha<\lambda$.
Let $B_\alpha$ be the set of all sets $b\in\roth$ such that
\begin{equation}\label{eq:h}
x\cap \Un_{n\in b}[a(n),a(n+1))\in\Fin,
\end{equation}
for all but less than $\kappa$ sets $x\in X_\alpha$.
The set $B_\alpha$ is groupwise dense:
Fix a set $c\in\roth$.
There is a set $b\in\roth$ such that 
\[
x\cap \Un_{n\in b}[(a\circ c)(n),(a\circ c)(n+1))\in\Fin,
\]
for all but less than $\kappa$ sets $x\in X$.
Let 
\[
d:=\Un_{n\in b}[c(n),c(n+1)).
\]
Fix a natural number $i\in d$.
There is a natural number $n\in b$ with $i\in[c(n),c(n+1))$.
Then $[a(i),a(i+1))\sub [a(c(n)), a(c(n+1)))$.
Thus, 
\[
\Un_{i\in d}[a(i),a(i+1))\sub \Un_{n\in b}[a(c(n)), a(c(n+1)))=\Un_{n\in b}[(a\circ c)(n),(a\circ c)(n+1)).
\]
It follows that $d\in B_\alpha$.

Since $\lambda<\fg$, there is a set $b\in\bigcap_{\alpha<\lambda}B_\alpha$.
For each ordinal number $\alpha<\lambda$, there is a set $S_\alpha\sub X_\alpha$ with $\smallcard{S_\alpha}<\kappa$ such that condition~\eqref{eq:h} holds for all sets $x\in X_\alpha\sm S_\alpha$.
Then condition~\eqref{eq:h} holds for all sets $x\in X\sm\Un_{\alpha<\lambda}S_\alpha$.
Since $\smallcard{\Un_{\alpha<\lambda}S_\alpha}<\kappa$, 
the set $X$ is a $\kappa$-generalized tower.

(2)
Proceed as in~(1) with the exception that $\smallcard{\Un_{\alpha<\lambda}S_\alpha}<(\kappa\cdot\lambda)^+$. 
\epf

\section{Products of sets satisfying $\sgg$}

\setcounter{equation}{0}

For a class $\cA$ of covers of spaces and a space $X$, let $\cA(X)$ be the family of all covers of $X$ from the class $\cA$.
Let $\cA$, $\cB$ be classes of covers of spaces.
A space $X$ satisfies $\sone(\cA,\cB)$ if for each sequence $\eseq{\cU}\in\cA(X)$, there are sets $\seleseq{U}$ such that $\sset{U_n}{n\in\bbN}\in\cB(X)$.
A \emph{Borel} cover of a space is a cover whose members are Borel subsets of the space.
Let $\Gab$ be the class of all Borel $\gamma$-covers of spaces.
For a property of spaces $\bfP$, let $\non(\bfP)$ be the minimal cardinality of a set of reals with no property $\bfP$.
We have~\cite[Theorem~4.7]{coc2},~\cite[Theorem~27(2)]{CBC}
\[
\non(\sbgg)=\non(
\sgg)=\fb.
\]
In the Laver model, all sets of reals satisfying $\sgg$ have cardinality strictly smaller than $\fb$~\cite[Theorem~3.6]{BBC}, and thus they are trivial; it implies that the properties $\sbgg$ and $\sgg$ are equivalent.
In order to construct \emph{nontrivial} examples of sets of reals with the above properties, additional set-theoretic assumptions are needed.
A set of reals of cardinality at least $\fb$ is a \emph{$\fb$-Sierpi\'nski set}, if its intersection with every Lebesgue-null set has cardinality less than $\fb$.
Every $\fb$-Sierpi\'nski set satisfies $\sbgg$ and assuming that $\cov(\cN)=\cof(\cN)=\fb$, such a set exists (definitions of cardinal characteristics of the continuum and relations between them can be found in the Blass survey~\cite{blass}).
For functions $a,b\in\roth$, we write $a\les b$ if the set $\sset{n}{b(n)<a(n)}$ is finite.
A set $A\sub\roth$ is \emph{bounded} if there is a function $b\in\roth$ such that $a\les b$ for all functions $a\in A$, denoted $A\les b$.
A subset of $\roth$ is \emph{dominating} if for every function $a\in\roth$, there is a function $b$ from the set such that $a\les b$.
Let $\fd$ be the minimal cardinality of a dominating set in $\roth$.
Assuming $\fp=\fb$ or $\fb<\fd$, a $\fb$-unbounded tower exists~\cite[Lemma 2.2]{BBC}.
Also if for an infinite ordinal number $\kappa$, there is a $\kappa$-generalized tower, then there is a $\fb$-generalized tower~\cite[Proposition 2.4]{BBC}.
Each $\fb$-unbounded tower set satisfies $\sone(\Ga,\Ga)$
(\cite[Theorem 6]{alpha_i}, \cite[Proposition 2.5]{BBC}), but not $\sbgg$~\cite[Corollary~2.10]{BBC}.

A subfamily of $\roth$ is \emph{open} if it is closed under almost subsets, and it is \emph{dense} if each set from $\roth$ has an almost subset in the family.
The \emph{density number} $\fh$ is the minimal cardinality of a collection of open dense families in $\roth$ with empty intersection. 
Let $\add(\sgg)$ be the minimal cardinality of a family of sets of reals satisfying $\sgg$ whose union does not satisfy $\sgg$.
By the result of Scheepers~\cite[Theorem~5]{SchCp} we have $\add(\sgg)\geq\fh$.
Miller and Tsaban proved that, assuming $\add(\sgg)=\fb$, each $\fb$-unbounded tower set satisfies $\sgg$ in all finite powers.
By the result of Miller, Tsaban and Zdomskyy, the product space of a $\fb$-unbounded tower set and a set satisfying $\sbgg$, satisfies $\sgg$~\cite[Theorem~7.1]{gamma}.
In the forthcoming Theorem~\ref{thm:tgg}, we generalize the above results.
In the proof we develop methods of Miller, Tsaban and Zdomskyy~\cite[Theorem~7.1]{gamma} and combine them with tools invented by Tsaban and the first named author~\cite[Lemma~5.1]{ST}.

\bthm\label{thm:tgg}
The product space of finitely many $\fb$-generalized tower sets and a set satisfying $\sbgg$, satisfies $\sgg$.
\ethm

In order to prove Theorem~\ref{thm:tgg}, we need the following notions and auxiliary results.
Let $\overline{\bbN}$ be the set $\bbN\cup\{\infty\}$ with the discrete topology and $n<\infty$ for all natural numbers $n$.
In the space $\overline{\bbN}^\bbN$, define relation $\les$, analogously as in $\roth$.
For a set $x\in\PN$, let $x\comp:=\bbN\sm x$.

\blem\label{lem:sbgg}
Let $m$ be a natural number, $Y$ be a set satisfying $\sbgg$ and $\sset{U_n}{n\in\bbN}\in \Ga(\Fin^m\x Y)$ be a family of open sets in $\PN^m\x Y$.
\be
\item 
There are Borel functions $f,g\colon Y\to \overline{\bbN}^\bbN$ such that for all points $x\in (\roth)^m$ and $y\in Y$, and all natural numbers $n$:
\[
\text{If }x_i\cap [f(y)(n),g(y)(n))=\emptyset\text{ for all natural numbers }i\leq m,\text{ then }(x,y)\in U_n.
\]
\item There is an increasing function $c\in\NN$ such that for each point $y\in Y$, we have 
\[
c(n)\leq f(y)(c(n+1))\leq g(y)(c(n+1))<c(n+2),
\]
for all but finitely many natural numbers $n$.
\ee
\elem

\bpf
(1)
Define a function $d\colon Y\to \NN$ as follows:
Fix a point $y\in Y$.
Let $d(y)(1)$ be the minimal natural number such that 
\[
\{\emptyset\}^m\x\{y\}\sub \bigcap_{k\geq d(y)(1)}U_{k}
\]
and, for each natural number $i$, let $d(y)(i+1)$ be the minimal natural number  greater than $d(y)(i)$ such that 
\[
\Pof([1,i+1))^m\x\{y\}\sub \bigcap_{k\geq d(y)(i+1)}U_{k}.
\]
The function $d$ is Borel:
We have
\[
\sset{y\in Y}{ d(y)(1)=1}=\medset{y\in Y}{\{\emptyset\}^m\x\{y\}\sub\bigcap_{k\geq 1}U_k},
\]
and
\[
\sset{y\in Y}{ d(y)(1)=j+1}=\medset{y\in Y}{\{\emptyset\}^m\x\{y\}\sub\bigcap_{k\geq j+1}U_k\sm U_j},
\]
for all natural numbers $j$; the above sets are Borel.
Fix natural numbers $i,j$, and assume that the set 
\[
\sset{y\in Y}{ d(y)(i)=j}
\]
is Borel.
Fix a point $y\in Y$ such that $d(y)(i+1)=j+1$.
We have
\[
\Pof([1,i+1))^m\x\{y\}\sub \bigcap_{k\geq j+1}U_{k}.
\]
 If $\Pof([1,i+1))^m\x\{y\}\sub U_j$, then $d(y)(i)=j$, otherwise there is an element $t\in\Pof([1,i+1))^m$ with $(t,y)\notin U_{j}$.
Thus,
\begin{align*}
\sset{y\in Y}{d(y)(i+1)=j+1}=&
\bigcap_{t\in\Pof([1,i+1))^m}\smallmedset{y\in Y}{(t,y)\in \bigcap_{k\geq j+1}U_k}\cap\\
&
\Bigl(\smallmedset{y\in Y}{d(y)(i)=j} \cup\Un_{t\in\Pof([1,i+1))^m}\set{y\in Y}{(t,y)\notin U_{j}}\Bigr).
\end{align*}
Then  the set $\sset{y\in Y}{ d(y)(i+1)=j+1}$ is Borel.

For a point $y\in Y$, let $f(y)\in\overline{\bbN}^\bbN$ be a function such that:
\[
f(y)(n):=
\begin{cases}
\infty,&\text{ for }n\in [1,d(y)(1)),\\
i,&\text{ for }n\in[d(y)(i),d(y)(i+1)), i\in\bbN.
\end{cases}
\]
For natural numbers $i,n$, we have
\[
\sset{y\in Y}{f(y)(n)=\infty}=\sset{y\in Y}{n<d(y)(1)}
\]
and
\[
\sset{y\in Y}{f(y)(n)=i}=\sset{y\in Y}{d(y)(i)\leq n < d(y)(i+1)}.
\]
Since the function $d$ is Borel, the function $f\colon Y\to \overline{\bbN}^\bbN$, defined above, is Borel as well.

For each natural number $n$, let $\eseq{U^{(n)}}$ be an increasing sequence of clopen sets in $\PN^m\x Y$ such that $U_n=\Un_kU^{(n)}_k$.
Fix a point $y\in Y$.
Define a function $g(y)\in\overline{\bbN}^\bbN$ in the following way:
Let $n$ be a natural number.
If $f(y)(n)=\infty$, then $g(y)(n):=\infty$.
Assume that $f(y)(n)\in\bbN$.
There is a natural number $j>f(y)(n)$ such that
\[
\Pof([f(y)(n), j)\comp)^m\x \{y\}\sub U^{(n)}_{j}:
\]
Let $k$ be the minimal natural number with 
\[
\Pof([1,f(y)(n)))^m\x \{y\}\sub U^{(n)}_{k}.
\] 
Let $j$ be the minimal natural number with $j>f(y)(n)$ and 
\[
\Pof([f(y)(n),j)\comp)^m\x\{y\}\sub U^{(n)}_k.
\]
Since the sets $\eseq{U^{(n)}}$ are ascending, we have
\[
\Pof([f(y)(n), j)\comp)^m\x \{y\}\sub U^{(n)}_{j}.
\]
Let $g(y)(n)>f(y)(n)$ be the minimal natural number with
\[
\Pof([f(y)(n), g(y)(n))\comp)^m\x \{y\}\sub U^{(n)}_{g(y)(n)}.
\]
The set $\Pof([f(y)(n), g(y)(n))\comp)^m$ is compact, and thus there is an open set $V^y_n$ in $Y$ such that 
\[
\Pof([f(y)(n), g(y)(n))\comp)^m\x V^y_n\sub U^{(n)}_{g(y)(n)}.
\]

The function $g$ is Borel:
Fix a natural number $n$.
We have $\sset{y\in Y}{g(y)(n)=\infty}=\sset{y\in Y}{f(y)(n)=\infty}$.
Since the function $f$ is Borel, the latter set is Borel.
For a natural number $k$, we have
\[
\sset{y\in Y}{g(y)(n)=k}=\Un_{i<k}\sset{y\in Y}{f(y)(n)=i\wedge g(y)(n)=k}.
\]
In order to show that such a set is Borel, it is enough to show that for natural numbers $i,k$ with $i<k$, the set $\sset{y\in Y}{f(y)(n)=i\wedge g(y)(n)=k}$ is Borel:
Fix natural numbers $i,k$ with $i<k$.
Let $y\in Y$ be a point such that $f(y)(n)=i$ and $g(y)(n)=k$.
For each point $y'\in V^y_n$ such that $f(y')(n)=i$, we have 
\[
\Pof([f(y')(n),g(y)(n)\comp)^m\x\{y'\}\sub \Pof([f(y)(n), g(y)(n))\comp)^m\x V^y_n\sub U^{(n)}_{g(y)(n)}.
\]
By the minimality of the number $g(y')(n)$, we have 
\beq\label{eq:ptwl}
g(y')(n)\leq g(y)(n).
\eeq
If $k=i+1$, then $g(y')(n)=g(y)(n)$, and thus 
the set $\sset{y\in Y}{f(y)(n)=i\wedge g(y)(n)=k}$ is Borel (in fact, it is even open).
Now assume that for a natural number $k> i+1$, the set $\sset{y\in Y}{f(y)(n)=i\wedge g(y)(n)<k}$ is Borel.
By the above inequality~\eqref{eq:ptwl}, we have 
\begin{gather*}
\sset{y\in Y}{f(y)(n)=
	i\wedge g(y)(n)=k}=\\
\Bigl(\sset{y\in Y}{f(y)(n)=i}\cap
\Un\sset{V^y_n}{g(y)(n)=k}\Bigr)\sm \sset{y\in Y}{f(y)(n)=i\wedge g(y)(n)<k}.
\end{gather*}
Since the function $f$ is Borel, the set $\sset{y\in Y}{f(y)(n)=i\wedge g(y)(n)=k}$ is Borel.

(2)
Let $d$ and $g$ be the functions from the proof of (1).
For each point $y\in Y$, we have $g(y)(n)=\infty$ for finitely many natural numbers $n$.
Since the functions $d$ and $g$ are Borel and the set $Y$ satisfies $\sbgg$, there is an increasing function $b\in\NN$ such that $\sset{d(y),g(y)}{y\in Y}\les b$~\cite[Theorem~1]{CBC}.
We may assume that $b(1)\neq 1$.
Let $c(1):=b(1)$ and $c(n+1):=b(c(n))$ for all natural numbers $n$.
For each point $y\in Y$, we have 
\[
c(n)\leq f(y)(c(n+1))\leq g(y)(c(n+1))<c(n+2)
\]
for all but finitely many natural numbers $n$:
Fix a point $y\in Y$.
There is a natural number $l$ such that 
\[
d(y)(n)\leq b(n)\text{ and } g(y)(n)\leq b(n),
\]
for all natural numbers $n\geq l$.
Let $n$ be a natural number with $n\geq l$.
Since 
\[
d(y)(c(n))\leq b(c(n))=c(n+1)
\]
and the function $f(y)$ is nondecreasing for arguments greater than or equal to $d(y)(1)$,
we have
\[
c(n)=f(y)(d(y)(c(n)))\leq f(y)(c(n+1))\leq g(y)(c(n+1))\leq b(c(n+1))=c(n+2).\qedhere
\]
\epf

\blem\label{lem:small_mistake_allowed}
Let $X$ be a space such that for every sequence $\eseq{\cU}\in\Ga(X)$,
there are sets $\seleseq{U}$ and a space $X'\sub X$ satisfying $\sgg$ such that $\sset{U_n}{n\in\bbN}\in \Ga(X\sm X')$.
Then the space $X$ satisfies $\sgg$.
\elem

\bpf
Let $\cU_n=\sset{U^{(n)}_m}{m\in\bbN}$ for all natural numbers $n$.
We may assume that for all natural numbers $n,m$, we have
\[
U^{(n+1)}_m\sub U^{(n)}_m.
\]
For each natural number $n$, we have $\sset{U^{(n)}_m}{m\geq n}\in\Ga(X)$.
Then  there are a function $f\in\NN$ and a space $X'\sub X$ such that $\sset{U^{(n)}_{f(n)}}{n\in\bbN}\in\Ga(X\sm X')$.
The range of the function $f$ is infinite, and thus there is a function $d\in\roth$ such that $(f\circ d)\in\roth$.
Let $f':=f\circ d$.
Since $U^{(d(n))}_{f'(n)}\sub U^{(n)}_{f'(n)}$, we have $\sset{U^{(n)}_{f'(n)}}{n\in\bbN}\in\Ga(X\sm X')$.

Since the space $X'$ satisfies $\sgg$, the set
\[
\set{a\in\roth}{\sset{U^{(n)}_{a(n)}}{n\in\bbN}\in \Ga(X')}
\]
is open and dense~\cite[Theorem~5]{SchCp}.
Thus, there is a function $a\in\roth$ with $a\sub f'$ such that $\sset{U^{(n)}_{a(n)}}{n\in\bbN}\in \Ga(X')$.
There is a function $c\in\roth$ such that $a=f'\circ c$.
Then $\sset{U^{(c(n))}_{a(n)}}{n\in \bbN}\in\Ga(X\sm X')$.
Since $U^{(c(n))}_{a(n)}\sub U^{(n)}_{a(n)}$, we have $\sset{U^{(n)}_{a(n)}}{n\in\bbN}\in\Ga(X)$.
\epf

\begin{proof}[{Proof of Theorem~\ref{thm:tgg}}]
Let $Y$ be a set satisfying $\sbgg$.

Let $X$ be a $\fb$-generalized tower and $\eseq{\cU}\in \Ga((X\cup\Fin)\x Y)$ be a sequence of families of open sets in $\PN\x Y$ such that $\cU_k=\sset{U^{(k)}_n}{n\in\bbN}$ for all natural numbers $k$.
For each natural number $k$, let $f_k,g_k, c_k$ be functions from Lemma~\ref{lem:sbgg} applied to the family $\cU_k$.
Let $h\colon Y\to\NN$ be a function such that $h(y)(k)$ is the minimal natural number $n$ such that for all natural numbers $j$ with $j\geq n$, we have 
\[
c_k(j)\le f_{{}k}(y)(c_k(j+1))<g_{{}k}(y)(c_k(j+1))<c_k(j+2).
\]
Since the functions $f_k$ and $g_k$ are Borel, so is the function $h$, and hence there is an increasing function $z\in\NN$ such that $\sset{h(y)}{y\in Y}\les z$.
We may assume that
\[
c_{k+1}(z(k+1))>c_k(z(k)+2)
\]
for all natural numbers $k$.
Since the set $X$ is a $\fb$-generalized tower, there are a set $b\in\roth$ and a set $\tX\sub X$ with $\smallcard{\tX}<\fb$ such that 
\[
x\cap \Un_{k\in b} [c_k(z(k)), c_k(z(k)+2))\in\Fin
\]
for all sets $x\in X\sm \tX$.
	
We have $\sset{U^{(k)}_{c_k(z(k+1))}}{k\in\bbN}\in \Ga ((X\sm S)\x Y)$:
Fix points $x\in \dtX$ and $y\in Y$.
There is a natural number $l$ such that for all natural numbers $k\geq l$, we have
\[
c_k(z(k))\le f_{{}k}(y)(c_k(z(k)+1))<g_{k}(y)(c_k(z(k)+1))<c_k(z(k)+2),
\]
\[
x\cap [c_k(z(k)), c_k(z(k)+2))=\emptyset,
\]
and thus
\[
x\cap [f_k(y)(z(k+1)), g_k(y)(z(k+1)))=\emptyset.
\]
By Lemma~\ref{lem:sbgg}(1), 
\[
(x,y)\in \bigcap_{k\geq l}U^{(k)}_{z(k+1)}.
\]

Since $\smallcard{\tX}<\fb$ and $\add(\sbgg)=\fb$~\cite[Corollary~2.4]{add}, the product space $(\tX\cup\Fin)\x Y$ satisfies $\sbgg$.
By Lemma~\ref{lem:small_mistake_allowed}, the product space $(X\cup\Fin)\x Y$satisfies $\sgg$.

Let $m$ be a natural number with $m>1$ and assume that the statement is true for $m-1$ $\fb$-generalized tower sets.
Let $X_1,\dotsc, X_{m}$ be $\fb$-generalized towers in $\roth$ and
\[
\eseq{\cU}\in \Ga((X_1\cup\Fin)\x\dotsb\x(X_{m}\cup\Fin)\x Y)
\]
be a sequence of families of open sets in $\PN^{m}\x Y$ such that $\cU_k=\sset{U^{(k)}_n}{n\in\bbN}$ for all natural numbers $k$.
We proceed analogously to the previous case.
For each natural number $k$, let $f_k,g_k, c_k$ be functions from Lemma~\ref{lem:sbgg} applied to the family $\cU_k$.
Let $h\colon Y\to\NN$ be a function such that $h(y)(k)$ is the minimal natural number $n$ such that for all natural numbers $j$ with $j\geq n$, we have 
\[
c_k(j)\le f_{{}k}(y)(c_k(j+1))<g_{{}k}(y)(c_k(j+1))<c_k(j+2).
\]
Since the functions $f_k$ and $g_k$ are Borel, so is the function $h$, and hence there is an increasing function $z\in\NN$ such that $\sset{h(y)}{y\in Y}\les z$.
We may assume that
\[
c_{k+1}(z(k+1))>c_k(z(k)+2)
\]
for all natural numbers $k$.
By Lemma~\ref{lem:uscale}(1), the set $X:=\Un_{i\leq m}X_i$ is a $\fb$-generalized tower.
Then  there are a set $b\in\roth$ and a set $\tX\sub X$ with $\smallcard{\tX}<\fb$ such that 
\[
x\cap \Un_{k\in b} [c_k(z(k)), c_k(z(k)+2))\in\Fin
\]
for all sets $x\in X\sm \tX$.
	
We have $\sset{U^{(k)}_{c_k(z(k+1))}}{k\in\bbN}\in \Ga ((X\sm S)^{m}\x Y)$:
Fix points $x\in (X\sm S)^{m}$ and $y\in Y$.
There is a natural number $l$ such that for all natural numbers $k\geq l$, we have
\[
c_k(z(k))\le f_{{}k}(y)(c_k(z(k)+1))<g_{k}(y)(c_k(z(k)+1))<c_k(z(k)+2),
\]
\[
x_i\cap [c_k(z(k)), c_k(z(k)+2))=\emptyset\text{ for all natural numbers }i\leq m,
\]
and thus
\[
x_i\cap [f_k(y)(z(k+1)), g_k(y)(z(k+1)))=\emptyset \text{ for all natural numbers }i\leq m.
\]
By Lemma~\ref{lem:sbgg}(1), 
\[
(x,y)\in \bigcap_{k\geq l}U^{(k)}_{z(k+1)}.
\]
Since $\smallcard{\tX}<\fb$ and $\add(\sbgg)=\fb$~\cite[Corollary~2.4]{add}, the product space $((X_j\cap\tX)\cup\Fin)\x Y$ satisfies $\sbgg$ for all natural numbers $j\leq m+1$, and the set
\[
\Un_{j\leq m}\prod_{\substack{i\leq m\\ i\neq j}}(X_i\cup \Fin)\x ((X_j\cap S)\cup\Fin)\x Y
\]
satisfies $\sgg$.
By Lemma~\ref{lem:small_mistake_allowed}, the product space $(X_1\cup\Fin)\x\dotsb\x(X_{m}\cup\Fin)\x Y$satisfies $\sgg$.
\epf

\bcor
Each $\fb$-generalized tower set satisfies $\sgg$ in all finite powers.
\ecor

The properties $\sgg$ and $\sbgg$ are closely related to local properties of functions spaces.
Let $X$ be a space and $\Cp(X)$ be the set of all continuous real-valued functions on $X$ with the pointwise convergence topology.
A sequence $\eseq{f}\in\Cp(X)$ converges \emph{quasinormally} to the constant zero function $\bfzero$, if there is a sequence of positive real numbers $\eseq{\epsilon}$ converging to zero such that for any point $x\in X$, we have $|f_n(x)|<\epsilon_n$ for all but finitely many natural numbers $n$.
A space $X$ is a \emph{QN-space} (\emph{wQN-space}) if every sequence $\eseq{f}\in\C(X)$ converging pointwise to $\bfzero$, converges (has a subsequence converging) quasinormally to $\bfzero$.
By a breakthrough result of Tsaban and Zdomskyy~\cite[Theorem~2]{TsZdArh},
a space is a QN-space if and only if all Borel images of the space in $\NN$ are bounded~\cite[Theorem 1]{CBC}.
In consequence, a space satisfies $\sbgg$ if and only if it is a QN-space.
Every perfectly normal space satisfying $\sgg$, is a wQN-space~\cite[Theorem~7]{wQN}.
QN-spaces, wQN-spaces and their variations were extensively studied by Bukovsk\'y, Hale\v{s}, Rec\l{}aw, Sakai and Scheepers~\cite{QN, wQN, brp, hales, sakaiCp, sakaisemcont, SchCp}.
We have the following corollary from Theorem~\ref{thm:tgg}.

\bcor
The product space of finitely many $\fb$-generalized tower sets and a QN-space, is a wQN-space.
\ecor

These results can also be formulated as dealing with Arhangel'ski\u{\i}'s properties
$\alpha_i$ for spaces of continuous real-valued functions~\cite{alpha_i, arch, arch2}.

\section{Generalized towers and productivity of $\Om\choose\Ga$}

\setcounter{equation}{0}

A space is \emph{Fr\'echet--Urysohn} if each point in the closure of a set is a limit of a sequence from the set.
By the celebrating result of Gerlits and Nagy, for a set of reals $X$,  the space $\Cp(X)$ is Fr\'echet--Urysohn if and only if the set $X$ satisfies $\GNga$~\cite[Theorem~2]{gn}.
The property $\GNga$ is \emph{productive} if each space satisfying $\GNga$, is productively $\GNga$.
The property $\GNga$ is preserved by finite powers, but
productivity of $\GNga$ is independent of ZFC.
In the Laver model, all sets of reals satisfying $\GNga$ are countable~(\cite[Theorem~2]{gn},~\cite[Theorem~17]{coc1},~\cite{Laver}), and thus $\GNga$ is productive; if \CH{} holds, then $\GNga$ is not productive~\cite[Theorem~3.2.]{gamma}.
Miller, Tsaban and Zdomskyy proved that, each $\w_1$-unbounded tower set is productively $\GNga$~\cite[Theorem~2.8.]{gamma}.
The main result of this section is the following Theorem.

\bthm\label{thm:mainGa}
\mbox{}
\be
\item
The product space of finitely many $\fp$-generalized tower sets, satisfies $\GNga$.
\item Assume that there is a $\fp$-generalized tower in $\roth$.
The following assertions are equivalent:
\be
\item Each set of reals of cardinality smaller than $\fp$ is productively $\GNga$.
\item Each $\fp$-generalized tower set is productively $\GNga$.
\ee
\ee
\ethm

Let $\Om$ be the class of all $\w$-covers of spaces.
 
\blem[{Galvin--Miller\cite[Lemma~1.2]{gami}}]\label{lem:GM}
Let $\cU$ be a family of open sets in $\PN$ such that $\cU\in\Om(\Fin)$.
There are a function $a\in\roth$ and sets $\eseq{U}\in \cU$ such that for each set 
$x\in \roth$ and all natural numbers~$n$:
\[
\text{If }x\cap[a(n), a(n+1))=\emptyset,\text{ then }x\in U_n.
\]
\elem

For spaces $X$ and $Y$, let $X\sqcup Y$ be the \emph{disjoint union} of these spaces.
The product space $X\x Y$ satisfies $\GNga$ if and only if $X\sqcup Y$ satisfies $\GNga$~\cite[Proposition~2.3]{miller}.
For functions $a,b\in\roth$, we write $a\leq b$ if $a(n)\leq b(n)$ for all natural numbers $n$.

\blem\label{lem:mainGa}
Let $X\sub\roth$ be a $\fp$-generalized tower and $Y$ be a set such that for every subset $\tX\sub X$ with $\smallcard{\tX}<\fp$, the product space $\tX\x Y$ satisfies $\GNga$. 
Then the product space $(X\cup\Fin)\x Y$ satisfies $\GNga$.
\elem

\bpf
Let  $\cU\in\Om((X\cup\Fin)\sqcup Y)$ be a family of open sets in $\PN\sqcup \PN$.
Let $\tX_1:=\Fin$.
Fix a natural number $k>1$, and assume that the set $\tX_{k-1}\sub X$  with $\Fin\sub \tX_{k-1}$ and $\smallcard{\tX_{k-1}}<\fp$ has been already defined.
Since $\smallcard{\tX_{k-1}}<\fp$, there is a subfamily $\cU'$ of $\cU$ with $\cU'\in\Ga(\tX_{k-1}\sqcup Y)$.
Apply Lemma~\ref{lem:GM} to the family $\cU'$.
Then there are a function $a_k\in\roth$ and sets $\eseq{U^{(k)}}\in\cV$ such that for each set 
$x\in \roth$ and all natural numbers $n$:
\beq\label{eq:GM}
\text{If }x\cap[a_k(n), a_k(n+1))=\emptyset,\text{ then }x\in U^{(k)}_n.
\eeq
Since the set $X$ is a $\fp$-generalized tower, there are a set $b_k\in\roth$ and a set $\tX_k\sub X$ with $\tX_{k-1}\sub \tX_k$ and $\smallcard{\tX_k}<\fp$ such that
\[
x\cap \Un_{n\in b_k}[a_k(n),a_k(n+1))\in\Fin
\]
for all sets $x\in X\sm \tX_k$.
Since $\cU'\in\Ga(Y)$, we have 
\[
\sset{U^{(k)}_{b_k(j)}}{j\in \bbN}\in\Gamma(((\dtX_k)\cup \tX_{k-1})\sqcup Y).
\]

There is a function $a\in\roth$ such that for each natural number $k$, we have 
\[
\card{(a_k\circ b_k)\cap [a(n),a(n+1))}\geq 2,
\]
for all but finitely many natural numbers $n$.
Since the set $X$ is a $\fp$-generalized tower, there are a set $b\in\roth$ and a set $\tX\sub X$ with $\smallcard{\tX}<\fp$ such that
\beq\label{eq:dtX}
x\cap \Un_{n\in b}[a(n),a(n+1))\in\Fin.
\eeq
for all sets $x\in X\sm \tX$.
We may assume that  $\Un_k \tX_k \sub \tX$.
The sets 
\beq\label{eq:b'}
b_k':=\medset{i\in b_k}{[a_k(i),a_k(i+1))\sub\Un_{n\in b} [a(n),a(n+1))}
\eeq
are infinite for all natural numbers $k$.
Thus,
\[
\set{U^{(k)}_{b_k'(j)}}{j\in\bbN}\in\Gamma(((\dtX_k)\cup \tX_{k-1})\sqcup Y).
\]
Since the sequence of the sets $\tX_k$ is increasing, we have $X=\Un_k(\dtX_k)\cup \tX_{k-1}$ and each point of $X$ belongs to all but finitely many sets $(\dtX_k)\cup \tX_{k-1}$.
For each point $x\in \tX$, define 
\[
g_x(k):=
\begin{cases}
 0, & \text{ if } x\nin (\dtX_k)\cup \tX_{k-1},\\
 \min\sset{j}{x\in\bigcap_{i\geq j}U^k_{b_k'(i)}}, & \text{ if } x\in (\dtX_k)\cup \tX_{k-1}.
\end{cases}
\]
Since $\smallcard{\tX}<\fp$, there is a function $g\in\NN$ with $\sset{g_x}{x\in \tX}\les g$ and
\beq\label{eq:g}
a_k(b_k'(g(k)+1))<a_{k+1}(b_{k+1}'(g(k+1))),
\eeq
for all natural numbers $k$.
Let 
\[
\cW_k:=\set{U^{(k)}_{b_k'(j)}}{j\geq g(k)}
\]
for all natural numbers $k$.
Then $\eseq{\cW}\in\Ga(Y)$.
We may assume that families $\cW_k$ are pairwise disjoint.
Since the properties $\GNga$ and $\swg$ are equivalent~\cite[Theorem~2]{gn}, the set $Y$ satisfies $\swg$.
Then there is a function $h\in\NN$ such that $g\leq h$ and
\[
\set{U^{(k)}_{b'_k(h(k))}}{k\in\bbN}\in\Ga(\tX\sqcup Y).
\]
Fix a set $x\in \dtX$.
By~\eqref{eq:b'}, for each natural number $k$, we have   
\[ 
\Un_{n\in b_k'}[a_k(n),a_k(n+1))\sub \Un_{n\in b}[a(n),a(n+1)).
\]
By~\eqref{eq:dtX}, \eqref{eq:g} and the fact that $g\leq h$, the set $x$ omits all but finitely many intervals
\[
\bigl[a_k(b_k'(h(k))),a_k(b_k'(h(k))+1)\bigr).
\]
By~\eqref{eq:GM}, we have 
\[
\set{U^{(k)}_{b'_k(h(k))}}{k\in\bbN}\in\Ga(\dtX).
\]
Thus, 
\[
\set{U^{(k)}_{b'_k(h(k))}}{k\in\bbN}\in\Ga((X\cup\Fin)\sqcup Y).\qedhere
\]
\epf

\begin{proof}[Proof of Theorem~\ref{thm:mainGa}]
(1) 
We prove a formally stronger assertion that the product space of finitely many $\fp$-generalized tower sets and a set of cardinality less than $\fp$, satisfies $\GNga$.
By~(a) and Lemma~\ref{lem:mainGa}, the product space of a $\fp$-generalized tower set and a set of cardinality smaller than $\fp$, satisfies $\GNga$.
Fix a natural number $m>1$.
Let $X_1,\dotsc, X_{m}$ be $\fp$-generalized towers in $\roth$ and $Y$ be a set with $\card{Y}<\fp$.
Assume that the product space 
\[
Z:=(X_1\cup\Fin)\x\dotsb\x(X_{m-1}\cup\Fin)\x Y
\]
satisfies $\GNga$.
Fix a set $\tX\sub X_{m}$ with $\smallcard{\tX}<\fp$. Since $\smallcard{\tX\x Y}<\fp$, by the inductive assumption, the product space
\[
\tX\x Z=(X_1\cup\Fin)\x\dotsb\x(X_{m-1}\cup\Fin)\x (\tX\x Y)
\]
satisfies $\GNga$.
By Lemma~\ref{lem:mainGa}, the product space $(X_1\cup\Fin)\x\dotsb\x(X_{m}\cup\Fin)\x Y$ satisfies $\GNga$.

(2)
($\Rightarrow$)
Apply Lemma~\ref{lem:mainGa}.

($\Leftarrow$)
Let $A\sub\roth$ be a set with $\card{A}<\fp$ and $Y$ be a set satisfying $\GNga$.
Since $\card{A}<\fp$, there is an element $b\in\roth$ such that $A\les b$.
Let $X\sub\roth$ be a $\fp$-generalized tower such that $x\sub b$ for all sets $x\in X$.
Then the set $Z:=X\cup A$ is a $\fp$-generalized tower.
We have 
\[
A=Z\cap\sset{x\in\roth}{b\les x},
\]
and thus the set $A$ is an $F_\sigma$ subset of $Z$.
The property $\GNga$ is preserved by taking $F_\sigma$ subsets~\cite[Theorem~3]{gami}.
Since the space $(Z\cup\Fin)\sqcup Y$ satisfies $\GNga$, the space $A\sqcup Y$ satisfies $\GNga$, too.
\epf

\section{Products of sets satisfying $\Om\choose\Ga$ and the properties $\sgg$ and $\sww$}

We already mentioned that the properties $\sgg$ and $\GNga$ were considered in the context of local properties of functions spaces. This is also the case for property $\sww$.
A space $Y$ has \emph{countable strong fan tightness}~\cite{sakai} if for each point $y\in Y$ and each sequence $\eseq{A}$ of subsets of the space $Y$ with $y\in\bigcap_n\overline{A_n}$, there are points $a_1\in A_1, a_2\in A_2,\dotsc$ such that $y\in\overline{\sset{a_n}{n\in\bbN}}$.
For a set of reals $X$, the space $\Cp(X)$ has countable strong fan tightness if and only if the set $X$ satisfies $\sww$~\cite{sakai}.

Let $\Op$ be the class of all open covers of spaces.
A space $X$ satisfies Menger's property $\men$ if for each sequence $\eseq{\cU}\in\Op(X)$, there are finite sets $\cF_1\sub\cU_1,\cF_2\sub \cU_2,\dotsc$ such that $\Un_n\cF_n\in\Op(X)$.
In this section, we consider products of sets satisfying $\GNga$ and their relations to the properties $\sgg$, $\sww$ and $\men$.
We have the following implications between considered properties.
\[
\xymatrix{
\sgg\ar[r] & \men\\
\GNga\ar[u]\ar[r]& \sww \ar[u]
}
\]
The properties $\sgg$ and $\sww$ are much stronger than $\men$.
Indeed, all sets of reals satisfying $\sgg$ or $\sww$ are \emph{totally imperfect}~\cite[Theorem~2.3]{coc2}, i.e., they do not contain uncountable compact subsets; each compact set satisfies $\men$.
The existence of a nontrivial set of reals satisfying $\sww$ is independent of $ZFC$: In the Laver model all sets of reals satisfying $\sww$ are countable(~\cite[Theorem~17]{coc1},~\cite{Laver}), and assuming that \CH{} holds, there is a nontrivial set satisfying $\sww$~\cite[Theorem~2.13]{coc2}.
In ZFC, there is a nontrivial totally imperfect set satisfying $\men$~\cite[Theorem~16]{BaTs}.

Miller, Tsaban and Zdomskyy proved that there are two sets of reals satisfying $\GNga$ whose product space does not satisfy $\men$~\cite[Theorem~3.2.]{gamma}.
In Theorem~\ref{thm:mainGa}, we gave a necessarily and sufficient condition when a $\fp$-generalized tower set is productively $\GNga$.
Now, we show that the product space of a $\fp$-generalized tower set and a set satisfying $\GNga$, satisfies $\sgg$ and $\sww$ in all finite powers (in fact, the property $\sww$ is preserved by finite powers~\cite{sakai}).

Let $\bfP$ be a property of spaces.
A space is \emph{productively $\bfP$}, if its product space with any space satisfying $\bfP$, satisfies $\bfP$.

\bthm\label{thm:mainsgg}
Let $\kappa$ be an uncountable ordinal number such that each set of reals of cardinality less than $\kappa$ is productively $\sgg$.
The product space of finitely many $\kappa$-generalized tower sets and a set satisfying $\GNga$, satisfies $\sgg$.
\ethm

Let $X$, $Y$ be spaces and $m$ be a natural number.
Identifying the space $X$ with $X\sqcup\emptyset$ and the space $Y$ with $\emptyset\sqcup Y$, the product space $X^m\x Y$ is a closed subset of $(X\sqcup Y)^{m+1}$.

\blem\label{lem:disjOm}
Let $X,Y\sub\PN$ and $m$ be a natural number.
Let $\cU\in\Om(\Fin^m\x Y)$ be a family of open sets in $\PN^{m+1}$.
There is a family $\cV\in\Om(\Fin\sqcup Y)$ of open sets in $\PN\sqcup\PN$ such that the family $\sset{V^{m+1}\cap (\PN^m\x Y)}{V\in\cV}$ refines the family $\cU$.
\elem

\bpf
There are families $\cW$ and $\cW'$ of open sets in $\PN$ such that the family 
\[
\sset{W^m\x W'}{W\in\cW, W'\in\cW'}\in \Om(\Fin^m\x Y)
\]
refines the family $\cU$.
We have 
\[
\sset{W\sqcup W'}{W\in\cW, W'\in\cW'}\in\Om(\Fin\sqcup Y),
\] 
and 
\[
W^m\x W'=(W\sqcup W')^{m+1}\cap (\PN^m\x Y)
\]
for all sets $W\in\cW$, $W'\in\cW'$.
\epf

\blem\label{lem:disjGa}
Let $\kappa$ be an uncountable ordinal number.
Let $X\sub\roth$ be a $\kappa$-generalized tower and $Y$ be a set satisfying $\GNga$.
For each family $\cU\in\Om(\Fin\sqcup Y)$ of open sets in $\PN\sqcup\PN$, there are a set $\tX\sub X$ with $\smallcard{\tX}<\kappa$ and a subfamily $\cV$ of the family $\cU$ with $\cV\in\Ga((X\sm\tX)\sqcup Y)$.
\elem

\bpf
Let $\cU\in\Om(\Fin\sqcup Y)$ be a family of open sets in $\PN\sqcup\PN$.
The space $\Fin\sqcup Y$ satisfies $\GNga$, and thus there is a subfamily $\cU'$ of the family $\cU$ with $\cU'\in\Ga(\Fin\sqcup Y)$.
Apply Lemma~\ref{lem:GM} to the family $\cU'$.
Then there are a function $a\in\roth$ and sets $\eseq{U}\in\cU'$ such that for each set 
$x\in \roth$ and all natural numbers $n$:
\[
\text{If }x\cap[a(n), a(n+1))=\emptyset,\text{ then }x\in U_n.
\]
Since $X$ is a $\kappa$-generalized tower, there are a set $b\in\roth$ and a set $\tX\sub X$ with $\smallcard{\tX}<\kappa$ such that 
\[
x\cap \Un_{n\in b}[a(n),a(n+1))\in\Fin
\]
for all sets $x\in X\sm \tX$.
We have $\sset{U_n}{n\in b}\in \Ga(\dtX)$.
Thus,
\[
\sset{U_n}{n\in b}\in\Ga((\dtX)\sqcup Y).\qedhere
\]
\epf

\blem\label{lem:partGa}
Let $\kappa$ be an uncountable ordinal number.
Let $X\sub\roth$ be a $\kappa$-generalized tower,  $Y$ be a set satisfying $\GNga$, and $m$ be a natural number.
For each sequence $\eseq{\cU}\in \Om(\Fin^m\x Y)$ of families of open sets in $\PN^{m+1}$, there are a set $\tX\sub X$ with $\smallcard{\tX}<\kappa$ and sets $\seleseq{U}$ such that  
\[
\sset{U_n}{n\in\bbN}\in\Gamma\bigl(((\dtX)\cup\Fin)^m\x Y\bigr).
\]
\elem

\bpf
Fix a natural number $m$. 
We may assume that the family $\cU_{n+1}$ refines the family $\cU_n$ for all natural numbers $n$.
By Lemma~\ref{lem:disjOm}, there is a sequence $\eseq{\cV}\in\Om(\Fin\sqcup Y)$ such that the family
\[
\sset{V^{m+1}\cap (X^m\x Y)}{V\in\cV_n}
\]
refines the family $\cU_n$ for all natural numbers $n$.
Since the set $\Fin\sqcup Y$ satisfies $\GNga$ and the properties $\GNga$ and $\swg$ are equivalent~\cite[Theorem~2]{gn}, there are sets $\seleseq{V}$ such that $\sset{V_n}{n\in\bbN}\in\Ga(\Fin\sqcup Y)$.
By Lemma~\ref{lem:disjGa}, there are a set $a\in\roth$ and a set $\tX\sub X$ with $\smallcard{\tX}<\kappa$ such that $\sset{V_n}{n\in a}\in\Ga(((\dtX)\cup\Fin)\sqcup Y)$.
For each natural number $n\in a$, there is a set $U_n\in\cU_n$ such that $V^{m+1}\cap (X^m\x Y)\sub U_n$.
For each natural number $n\in a\comp$, there is a set $U_n\in\cU_n$ such that $U_n\supseteq U_k$ for some natural number $k\in a$ with $n<k$.
We have $\sset{U_n}{n\in\bbN}\in \Ga(((\dtX)\cup\Fin)^m\x Y)$.
\epf

\begin{proof}[{Proof of Theorem~\ref{thm:mainsgg}}]
Let $Y$ be a set satisfying $\GNga$.

Let $X\sub\roth$ be a $\kappa$-generalized and $\eseq{\cU}\in \Ga((X\cup\Fin)\x Y)$ be families of open sets in $\PN^2$. 
By Lemma~\ref{lem:partGa}, there are sets $\seleseq{U}$ and a set $\tX\sub X$ with $\smallcard{\tX}<\kappa$ such that
\[
\set{U_n}{n\in\bbN}\in\Ga(((\dtX)\cup\Fin)\x Y).
\]
Since $\smallcard{\tX}<\kappa$ and, by the assumption, each set of cardinality smaller than $\kappa$ is productively $\sgg$, the product space $(\tX\cup\Fin)\x Y$ satisfies $\sgg$.
By Lemma~\ref{lem:small_mistake_allowed}, the product space $(X\cup\Fin)\x Y$ satisfies $\sgg$.

Fix a natural number $m>1$ and assume that the statement is true for $m-1$ $\kappa$-generalized tower sets.
Let $X_1,\dotsc,X_m\sub\roth$ be $\kappa$-generalized tower sets and 
\[
\eseq{\cU}\in \Ga((X_1\cup\Fin)\x\dotsb\x (X_m\cup\Fin)\x Y)
\]
be families of open sets in $\PN^{m+1}$. 
By Lemma~\ref{lem:uscale}(1), the set $X:=\Un_{i\leq m}X_i$ is a $\kappa$-generalized tower.
By Lemma~\ref{lem:partGa}, there are sets $\seleseq{U}$ and a set $\tX\sub X$ with $\smallcard{\tX}<\kappa$ such that
\[
\set{U_n}{n\in\bbN}\in\Ga(((\dtX)\cup\Fin)^m\x Y).
\]

The set
\[
Z:=
((X_1\cup\Fin)\x\dotsb\x (X_m\cup\Fin)\x Y)\sm (((\dtX)\cup\Fin)^m\x Y)
\]
satisfies $\sgg$:
Fix a natural number $i\leq m$.
By the inductive assumption, the product space
\[
\prod_{\substack{j\leq m\\ j\neq i}} (X_j\cup\Fin)\x Y
\]
satisfies $\sgg$.
By the assumption, each set of cardinality smaller than $\kappa$ is productively $\sgg$.
Since $\smallcard{X_i\cap \tX}<\kappa$,
the product space
\[
\prod_{\substack{j\leq m\\ j\neq i}} (X_j\cup\Fin)\x (X_i\cap \tX) \x Y
\]
satisfies $\sgg$.
A finite union of spaces satisfying $\sgg$, satisfies $\sgg$~\cite[Theorem~5]{SchCp}, and thus the set
\[
Z=\Un_{i\leq m}\prod_{\substack{j\leq m\\ j\neq i}} (X_j\cup\Fin)\x (X_i\cap \tX) \x Y,
\]
satisfies $\sgg$, too.

By Lemma~\ref{lem:small_mistake_allowed}, the product space $(X_1\cup\Fin)\x\dotsb\x (X_m\cup\Fin)\x Y$ satisfies $\sgg$.
\epf

The property $\GNga$ is preserved by finite powers and each set of cardinality less than $\fp$ is productively $\sgg$~\cite[Proposition~6.8]{blass},~\cite[Theorem~5]{SchCp}.
Thus, we have the following result.

\bcor
The product space of a $\fp$-generalized tower set and a set satisfying $\GNga$, satisfies $\sgg$ in all finite powers.
\ecor

A set $X\sub\NN$ is \emph{guessable} if there is a function $a\in\NN$ such that the sets $\set{n}{a(n)=x(n)}$ are infinite, for all functions $x\in X$.
Let $\cov(\cM)$ be the minimal cardinality of a family of meager subsets of the Baire space $\NN$, that covers $\NN$.
The minimal cardinality of a subset of $\NN$ that is no guessable, is equal to $\cov(\cM)$~\cite[Theorem~5.9]{blass}.

\bthm\label{thm:sww}
Let $\kappa$ be an uncountable ordinal number with $\kappa\leq\cov(\cM)$ such that $\kappa$ is regular or $\kappa<\cov(\cM)$.
The product space of finitely many $\kappa$-generalized tower sets and a set satisfying $\GNga$, satisfies $\sww$. 
\ethm

We need the following Lemma.

\blem\label{lem:sumaY}
A union of less than $\cov(\cM)$ sets satisfying  $\GNga$, satisfies  $\swo$.
\elem

\bpf
Fix an ordinal number $\lambda<\cov(\cM)$.
Let $X:=\bigcup_{\alpha<\lambda} X_\alpha$ be a union of sets satisfying  $\GNga$ and $\eseq{\cU}\in\Om(X)$. 
Assume that $\cU_n=\sset{U_m^{(n)}}{m\in\bbN}$ for all natural numbers $n$.
For each ordinal number $\alpha<\lambda$, there is a function $f_\alpha\in\NN$ such that $\sset{U_{f_\alpha(n)}^n}{n\in\bbN}\in\Ga(X_\alpha)$.
Since $\lambda<\cov(\cM)$, the set $\sset{f_\alpha}{\alpha<\kappa}$ is guessable.
Then there is a function $g\in\NN$ such that the sets $\sset{n}{f_\alpha(n)=g(n)}$ are infinite, for all ordinal numbers $\alpha<\kappa$.
Thus, $\sset{U^{(n)}_{g(n)}}{n\in\bbN}\in\Op(X)$.
\epf

\bpf[Proof of Theorem~\ref{thm:sww}]
A set of reals satisfies $\sww$ if and only if it satisfies $\soo$ in all finite powers~\cite{sakai}.
The properties $\soo$ and $\swo$ are equivalent~\cite[Theorem~17]{coc1} and property $\GNga$ is preserved by finite powers.
Thus, it is enough to show that the statement is true, when consider the property $\swo$ instead of $\sww$.
We prove a formally stronger assertion that the product space of finitely many $\kappa$-generalized tower sets, a set satisfying $\GNga$ and a set of cardinality less than $\cov(\cM)$, satisfies $\swo$.

Let $Y$ be a set satisfying $\GNga$ and $Z$ be a set with $\card{Z}<\cov(\cM)$.

Let $X\sub\roth$ be a $\kappa$-generalized tower and $\eseq{\cU}\in\Om((X\cup\Fin)\x Y\x Z)$ be families of open sets in $\PN^3$, where $\cU_n=\sset{U^{(n)}_m}{m\in\bbN}$ for all natural numbers $n$.
Fix a point $z\in Z$.
By Lemma~\ref{lem:partGa}, there are a function $g_z\in\NN$ and a set $S_z\sub X$ with $\card{S_z}<\kappa$ such that 
\[
\sset{U^{(n)}_{g_z(n)}}{n\in\bbN}\in \Ga(((X\sm S_z)\cup \Fin)\x Y \x \{z\})
\]
Since $\card{Z}<\cov(\cM)$, there is a function $g\in\NN$ such that the sets $\sset{n}{g(n)=g_z(n)}$ are infinite for all points $z\in Z$.
Let $S:=\Un\sset{S_z}{z\in Z}$.
We have
\[
\sset{U^{(n)}_{g(n)}}{n\in\bbN}\in \Ga(((X\sm S)\cup \Fin)\x Y \x Z)
\]
Then  there are sets $U_1\in\cU_1, U_3\in\cU_3,\dotsc$ such that 
\[
\sset{U_{2n-1}}{n\in
\bbN}\in \Ga(((X\sm S)\cup \Fin)\x Y \x Z).
\]
By the assumption about the ordinal number $\kappa$, we have $\card{S}<\cov(\cM)$.
Thus, $\card{S\x Z}<\cov(\cM)$.
By Lemma~\ref{lem:sumaY}, the product space $S\x Y \x Z$ satisfies $\swo$.
There are sets $U_2\in\cU_2, U_4\in\cU_4,\dotsc$ such that 
\[
\sset{U_{2n}}{n\in\bbN}\in\Op(S\x Y\x Z).
\]
Finally, we have
\[
\sset{U_{n}}{n\in\bbN}\in\Op((X\cup\Fin)\x Y\x Z).
\]

Fix a natural number $m>1$ and assume that the statement is true for $m-1$ $\kappa$-generalized tower sets.
Let $X_1,\dotsc, X_m\sub\roth$ be $\kappa$-generalized towers.
Let
\[
\eseq{\cU}\in\Om((X_1\cup\Fin)\x(X_m\cup\Fin)\x Y \x Z)
\]
be a sequence of families of open sets in $\PN^{m+2}$, where $\cU_n=\sset{U^{(n)}_m}{m\in\bbN}$ for all natural numbers $n$.
By Lemma~\ref{lem:uscale}(1), the set $X:=\Un_{i\leq m} X_i$ is a $\kappa$-generalized tower.
Fix a point $z\in Z$.
By Lemma~\ref{lem:partGa}, there are a function $g_z\in\NN$ and a set $S_z\sub X$ with $\card{S_z}<\kappa$ such that 
\[
\sset{U^{(n)}_{g_z(n)}}{n\in\bbN}\in \Ga(((X\sm S_z)\cup \Fin)^m\x Y \x \{z\}).
\]
Since $\card{Z}<\cov(\cM)$, there is a function $g\in\NN$ such that the sets $\sset{n}{g(n)=g_z(n)}$ are infinite for all points $z\in Z$.
Let $S:=\Un\sset{S_z}{z\in Z}$.
We have
\[
\sset{U^{(n)}_{g(n)}}{n\in\bbN}\in \Ga(((X\sm S)\cup \Fin)^m\x Y \x Z).
\]
Then  there are sets $U_1\in\cU_1, U_3\in\cU_3,\dotsc$ such that 
\[
\sset{U_{2n-1}}{n\in
\bbN}\in \Ga(((X\sm S)\cup \Fin)^m\x Y \x Z).
\]

The set
\[
T:=
((X_1\cup\Fin)\x\dotsb\x (X_m\cup\Fin)\x Y\x Z)\sm (((\dtX)\cup\Fin)^m\x Y\x Z)
\]
satisfies $\swo$:
Fix a natural number $i\leq m$.
By the inductive assumption, the product space
\[
\prod_{\substack{j\leq m\\ j\neq i}} (X_j\cup\Fin)\x Y\x Z
\]
satisfies $\swo$.
By the assumption about the ordinal number $\kappa$, we have $\card{S}<\cov(\cM)$.
Since $\smallcard{X_i\cap \tX}<\cov(\cM)$, we have $\smallcard{(X_i\cap \tX)\x Z}<\cov(\cM)$.
By the inductive assumption, the product space
\[
\prod_{\substack{j\leq m\\ j\neq i}} (X_j\cup\Fin)\x Y\x (X_i\cap \tX)\x Z
\]
satisfies $\swo$.
A finite union of sets satisfying $\swo$, satisfies $\swo$~\cite[Theorem~2.3.9]{BarJu}, and thus the set
\[
T=\Un_{i\leq m}\prod_{\substack{j\leq m\\ j\neq i}} (X_j\cup\Fin)\x Y\x (X_i\cap \tX)\x Z,
\]
satisfies $\swo$, too.
There are sets $U_2\in\cU_2, U_4\in\cU_4,\dotsc$ such that 
\[
\sset{U_{2n}}{n\in\bbN}\in\Op(T).
\]
Finally, we have
\[
\sset{U_{n}}{n\in\bbN}\in\Op((X_1\cup\Fin)\x\dotsb\x (X_m\cup\Fin)\x Y\x Z).\qedhere
\]
\epf

\bcor
Let $\kappa$ be an uncountable ordinal number with $\kappa\leq\cov(\cM)$ such that $\kappa$ is regular or $\kappa<\cov(\cM)$.
Each $\kappa$-generalized tower set satisfies $\sww$.
\ecor

Since the ordinal number $\fp$ is regular and $\fp\leq\cov(\cM)$, we have the following result.

\bcor
The product space of a $\fp$-generalized tower set and a set satisfying $\GNga$, satisfies $\sww$.
\ecor

\section{Remarks and open problems}

\subsection{Around Scheepers's Conjecture}

A \emph{clopen} cover of a space is a cover whose members are clopen subsets of the space.
Let $\Gaclp$ be the class of all clopen $\gamma$-covers of spaces.
A set of reals is a wQN-space if and only if it satisfies $\sclpgg$~\cite[Theorem~9]{QN}.
The following conjecture was formulated by Scheepers.

\bcnj[{\cite{SakaiScheepersPIT}}]\label{cnj:sch}
The properties $\sgg$ and $\sclpgg$ are equivalent.
\ecnj

The property $\sgg$ describes a local property of functions spaces:
Let $\bbR$ be the real line with the usual topology.
Let $X$ be a set of reals.
A function $f\colon X\to \bbR$ is \emph{upper semicontinuous} if the sets $\sset{x\in X}{f(x)<a}$ are open for all real numbers $a$.
By the result of Bukovsk\'y~\cite[Theorem~13]{wQN}, the set $X$ satisfies $\sgg$ if and only if it is an \emph{$\text{SPP}^*$ space}, that is, for each sequence $\seq{f_{1,m}}{m}, \seq{f_{2,m}}{m},\dotsc$ of sequences of upper continuous functions on $X$, each of them converging pointwise to the constant zero function $\bfzero$, there is a sequence $\seq{m_n}{n}$ of natural numbers such that the sequence $\seq{f_{n,m_n}}{n}$ converges to $\bfzero$.
Thus, in the language of functions spaces, the above Scheepers conjecture asks whether the properties wQN and $\text{SSP}^*$ are equivalent.

\subsection{Generalized towers and productivity of $\hur$}

A space $X$ satisfies Hurewicz's property $\hur$ if for each sequence $\eseq{\cU}\in\Op(X)$ there are finite sets $\cF_1\sub\cU_1, \cF_2\sub\cU_2,\dotsc$ such that $\sset{\Un\cF_n}{n\in\bbN}\in \Ga(X)$.
We have the following implications between considered properties
\[
\sgg\longrightarrow \hur \longrightarrow \men,
\]
and the property $\hur$ is strictly in between properties $\sgg$ and $\men$ (\cite[Theorems~2.2, 2.4]{coc2},~\cite[Theorem~3.9]{sfh}).
Let $\cFin$ be the set of all cofinite subsets of $\bbN$.
For functions $z,t\in \roth$, we write $z\leinf t$ if $t\not\les z$.

\bdfn[{\cite[Definition~4.1]{ST}}]
A set $X\sub\roth$ is a \emph{$\cFin$-scale} if for each element $z\in\roth$, there is an element $t\in\roth$ such that
\[
z\leinf t\leqslant^* x
\]
for all but less than $\fb$ functions $x\in X$.
\edfn

\blem
Each $\fb$-generalized tower is a $\cFin$-scale.
\elem

\bpf
Let $X\sub\roth$ be a $\fb$-generalized tower and $z\in\roth$ be an element such that $z(1)\neq1$.
Define an element $\tz$ such that $\tz(1):=z(1)$, and
$\tz(n+1):=z(\tz(n))$ for all natural numbers $n$.
There is a set $b\in\roth$ such that
\[
x\cap\Un_{n\in b} [\tz(n),\tz(n+1))\in\Fin
\]
for all but less than $\fb$ many elements $x\in X$.
We have $b\comp\in\roth$.
Then the set $t:=\bigcup_{n\in b\comp} [\tilde{z}(n),\tilde{z}(n+1))$ omits infinitely many intervals $[\tz(n),\tz(n+1))$ and $x\as t$ for all but less than $\fb$ elements $x\in X$.
We have 
\[
z(\tz(n))\leq\tz(n+1)\leq t(\tz(n)),
\]
for all natural numbers $n\in b\comp$.
Thus, $z\leqslant^\infty t$.
For an element $x\in\roth$ such $x\as t$ and $t\sm x\in\roth$, we have $t\les x$.
There are only countably many elements $x\in\roth$ with $t\sm x\in \Fin$.
Thus,
\[
z\leinf t\les x
\]
for all but less than $\fb$ elements $x\in X$.
\epf

By the result of Tsaban and the first named author~\cite[Theorem~5.4]{ST}, we have the following corollary.

\bcor
Let $\kappa$ be an uncountable ordinal number with $\kappa\leq \fb$.
Each $\kappa$-generalized tower set is productively $\hur$.
\ecor

\subsection{Questions}

\bprb
Is a $\fb$-unbounded tower, provably, productively $\sone(\Ga,\Ga)$?
Is this the case assuming \CH{}?
\eprb

\bprb
Assume Martin Axiom and the negation of \CH{}.
Is each set of cardinality less than $\fc$ productively $\GNga$?
\eprb

\bprb
Is it consistent that $\fp>\w_1$, each set of cardinality less than $\fp$ is productively $\GNga$ and there is a $\fp$-generalized tower?
\eprb

\bprb
Let $\kappa$ be an uncountable ordinal number.
Does the existence of a $\kappa$-generalized tower imply the existence of a $\kappa$-unbounded tower?
\eprb

\bprb
Let $\kappa$ be an uncountable ordinal number.
Is a union of less than $\fb$ many $\kappa$-generalized towers, a $\kappa$-generalized tower?
\eprb

\end{document}